\documentclass[%
 aip,
 amsmath,amssymb,
 reprint,%
]{revtex4-1}

\usepackage{graphicx}
\usepackage{dcolumn}
\usepackage{bm}

\usepackage[utf8]{inputenc}
\usepackage[T1]{fontenc}
\usepackage{mathptmx}
\usepackage{etoolbox}
\usepackage{xcolor}

\makeatletter
\def\@email#1#2{%
 \endgroup
 \patchcmd{\titleblock@produce}
  {\frontmatter@RRAPformat}
  {\frontmatter@RRAPformat{\produce@RRAP{*#1\href{mailto:#2}{#2}}}\frontmatter@RRAPformat}
  {}{}
}%
\makeatother

\newcommand{\di}[2]{\frac{d #1}{d #2}}

\begin{document}

\preprint{AIP/123-QED}

\title[Rate-induced biosphere collapse in the Daisyworld model]{Rate-induced biosphere collapse in the Daisyworld model}
\author{Constantin W. Arnscheidt}
 \affiliation{Centre for the Study of Existential Risk, University of Cambridge, Cambridge, UK, CB2 1SB}
  \email{ca628@cam.ac.uk}
\author{Hassan Alkhayuon}%
\affiliation{ 
School of Mathematical Sciences, University College Cork, Cork, Ireland, T12 XF62}

\date{\today}

\begin{abstract}
There is much interest in the phenomenon of rate-induced tipping, where a system changes abruptly when forcings change faster than some critical rate. Here, we demonstrate and analyse rate-induced tipping in the classic "Daisyworld" model. The Daisyworld model considers a hypothetical planet inhabited only by two species of daisies with different reflectivities, and is notable because the daisies lead to an emergent "regulation" of the planet's temperature. The model serves as a useful thought experiment regarding the co-evolution of life and the global environment, and has been widely used in the teaching of Earth system science.

We show that sufficiently fast changes in insolation (i.e. incoming sunlight) can cause life on Daisyworld to go extinct, even if life could in principle survive at any fixed insolation value among those encountered. 
Mathematically, this occurs due to the fact that the solution of the forced (nonautonomous) system crosses the stable manifold of a saddle point for the frozen (autonomous) system. 
The new discovery of rate-induced tipping in such a classic, simple, and well-studied model 
provides further supporting evidence that rate-induced tipping --- and indeed, rate-induced collapse --- may be common in a wide range of systems. 

\end{abstract}

\maketitle

\begin{quotation}
The concept of "tipping points" captures the idea that small changes can make a big difference. Rate-induced tipping is a kind of tipping which occurs when forcings (i.e. changing external conditions) change faster than a critical rate. In this article, we show that the classic "Daisyworld" dynamical-system model, which is essentially a thought experiment about co-evolution between life and the environment, exhibits rate-induced tipping towards collapse. This result is of both practical and fundamental interest: it makes the Daisyworld model more useful as a paedagogical tool, suggests that rate-induced tipping might be found in other classic models if we only look for it, and helps support the case that rate-induced tipping is a common phenomenon. 
\end{quotation}

\section{\label{sec:Introduction}Introduction}
When a system is pushed beyond a tipping point, change becomes self-perpetuating, and the state of the system may change dramatically. Tipping points occur in a wide range of systems, from environmental to societal \cite{folke04,lenton13_rev,scheffer20}, and have been the cause of much policy concern, especially in the context of anthropogenic climate change \cite{armstrongmckay22,gtp23}. While emphasis is often on tipping points in "natural" systems, they can also exist in human systems (societies, economies, infrastructure, etc.), with both potentially positive\cite{lenton22_postp,otto20,winkelmann22} and negative\cite{juhola22,spaiser24,arnscheidt24_sr} consequences. 

Past work on tipping points has often focused on those with fixed thresholds: in other words, external forcing needs to exceed some critical value before tipping occurs. Mathematically, this can occur due to "dangerous" bifurcations\cite{thompson1994safe} in the system, for example where a previously stable equilibrium state of the system loses its stability. However, such bifurcation-induced tipping, or "B-tipping", is only one way in which a system may tip\cite{ashwin12}.

Another kind of tipping, which has been the subject of much recent interest, is rate-induced tipping, or "R-tipping". Here, the system tips when external forcings change faster than a critical rate \cite{ashwin12,ritchie23}. Rate-induced tipping has been shown to be relevant in a wide range of contexts, from biology to climate to power grids \cite{stocker1997influence, scheffer08,wieczorek10, luke2011soil, siteur2016ecosystems, alkhayuon2019basin, okeeffe20, vanselow19,lohmann21,arnscheidt20,arnscheidt22_review,arnscheidt22_rate,ritchie23, feudel2023rate, panahi2023rate}. It occurs in simple dynamical system models, high-complexity ocean models\cite{lohmann21}, and agent-based ecological models\cite{arnscheidt22_rate}. It has been argued to be ubiquitous in a wide range of systems\cite{arnscheidt22_rate,ritchie23}. Indeed, given that the human-dominated Anthropocene is notable in part for its rapid rates of change\cite{steffen15}, rate-induced tipping may be a key concept for understanding the kinds of dangers humanity may face in the future. 

"Daisyworld" is a simple mathematical model first introduced by \textcite{watson83} as a thought experiment regarding the co-evolution of life and the global environment. It describes the dynamics of a hypothetical planet inhabited by two species of daisies: white and black. Remarkably, the system exhibits emergent self-regulation of temperature: the daisies essentially contribute to making the planet more habitable for themselves. In the four decades since it was introduced, the model has been the subject of much discussion and analysis, and has also spawned many extensions, variants, and applications; much of this work is reviewed by \textcite{wood08}. It has also often been used in the teaching of Earth system science\cite{mcguffie05_book,kump10_book,ford10_book}.

While the Daisyworld model and its offshoots stand independently as an interesting domain of study, we note that the original model was introduced in the context of Lovelock's Gaia hypothesis \cite{lovelock72,lovelock74} of planetary-scale homeostasis "by and for the biosphere". Here, Daisyworld is key because it shows that smaller-scale interactions between life and its environment can, in principle, lead to emergent global environmental regulation. 
Of course, this does not mean that such behaviour necessarily occurs in the real world. Some key arguments in the long-standing and ongoing debate on these issues are that Darwinian evolution provides no reason to expect such self-regulation to arise \cite{doolittle81,kirchner89,kirchner02}, but also more recently that there may be relevant selection mechanisms beyond Darwinian evolution \cite{doolittle14,lenton18selection,tamre24}. In this article we take no position on any of these issues, and will not discuss the Gaia hypothesis further. We will also not discuss the many extensions and variants of the Daisyworld model. 

Instead, we show for the first time that the original Daisyworld model of \textcite{watson83} exhibits rate-induced tipping. Specifically, Daisyworld's biosphere can collapse entirely (both species of daisies go extinct) if insolation (incoming sunlight) changes too quickly. This is the case even if the biosphere could in principle survive at any fixed insolation value among those encountered. The rate-induced tipping threshold emerges due to a "basin instability" \cite{okeeffe20}, in which the system state crosses the stable manifold of a saddle point. 

This result is of interest for a number of reasons. First, the fact that rate-induced tipping can be found in a model as well-known and well-studied as the Daisyworld model, four decades since its inception, suggests that rate-induced tipping may also be found in a wide range of other existing models if we only look for it. Second, since the emergent self-regulation in the Daisyworld model is a classic (simple) instance of co-evolution between life and its environment, this result helps support earlier arguments that rate-induced tipping (and even more specifically, rate-induced collapse) may be common in such contexts\cite{arnscheidt22_rate,arnscheidt22_review}. Finally, this increases the paedagogical value of the Daisyworld model. Daisyworld is already used to help teach Earth system science \cite{mcguffie05_book,kump10_book,ford10_book}, and could now be used to introduce students to an even richer range of dynamics, like rate-induced tipping, which are relevant for understanding various real-world systems.

\section{The Daisyworld model}
The Daisyworld model is specified by the two coupled ordinary differential equations
\begin{equation}
\di{\alpha_w}{t} = \alpha_w (\alpha_g \beta(T_w) - \gamma),
\label{eq:daw_dt}
·\end{equation}
\begin{equation}
\di{\alpha_b}{t} = \alpha_b (\alpha_g \beta(T_b) - \gamma).
\label{eq:dab_dt}
\end{equation}
Here $\alpha_w,\alpha_w$ are the fraction of the planet's surface covered by white and black daisies, respectively. $\alpha_g=1-\alpha_b-\alpha_w$ is the fraction of uncovered ground, while $\beta(T)$ is the growth rate and $\gamma$ the death rate of the daisies. It is assumed that there are no nutrient limitations, and that the entire surface of the planet is in principle suitable for hosting daisies. The growth rates of both daisy populations depend quadratically on temperature, with some optimum temperature $T_{\rm opt}$.
\begin{equation}
\beta(T) = \begin{cases}
1-k(T_{\rm opt} - T)^2, & |T-T_{\rm opt}|<1/\sqrt{k},\\
0, & \text{otherwise.} \\
\end{cases}
\label{eq:growth_rates}
\end{equation}
Here, $k$ is a parameter which determines how sensitive the growth rate is to temperature.

The effective emission temperature of the planet, $T_e$, is given by the radiative balance
\begin{equation}
\sigma T_e^4 = SL(1-A).
\label{eq:average_heat_balance}
\end{equation} 
$\sigma T_e^4$ is the radiation emitted by the planet, and $SL(1-A)$ is the radiation received from Daisyworld's host star.  $S$ is a radiative flux in units of W m$^{-2}$. We follow \textcite{watson83} in referring to $L$ as a dimensionless luminosity: scaling the star's luminosity by some factor is equivalent to scaling insolation (the actual radiation arriving at the planet's surface) by the same factor. 

For simplicity the planet is assumed to have no greenhouse gases in its atmosphere, and is assumed to be flat, such that the factor of $1/4$ that arises in energy balance equations for rapidly rotating spherical planets is neglected. $A$ is the albedo (reflectivity) of the planet: it is given by 
\begin{equation}
A = \alpha_w A_w + \alpha_b A_b + \alpha_g A_g,
\end{equation}
where $A_w$ and $A_b$ are the albedos of the white and black daisies and $A_g$ is the albedo of uncovered ground.
Finally, it is assumed that, due to the different albedos, each of the three kinds of surface (white daisy, black daisy, uncovered) has its own local temperature
\begin{equation}
T_i^4 = q(A-A_i)+T_e^4,
\label{eq:horiz_heat_balance}
\end{equation}
with $i$ replaced by $b,w$, or $g$ as appropriate.
$q$ is a positive constant quantifying the horizontal heat transfer across regions with different albedos. The outgoing radiative fluxes can thus be summed to recover the total outgoing flux.
\begin{equation}
	\alpha_w \sigma T_w^4 + \alpha_b \sigma T_b^4 + \alpha_g \sigma T_g^4 = \sigma T_e^4.
\end{equation}

To facilitate comparisons with previous studies, our parameter values are all the same as those of \textcite{watson83}, except for the heat transfer coefficient $q$ (see below). We use $\gamma=0.3$ (unit time)$^{-1}$, $k = 0.003265$ (unit time)$^{-1}$ K$^{-2}$, $T_{\rm opt} = 295.5$ K, $S=917$ W m$^{-2}$, $A_w=0.75$, $A_g=0.5$, and $A_b=0.25$. $\sigma$ is the Stefan-Boltzmann constant, which we take as 5.67 $\times$ 10$^{-8}$ W m$^{-2}$ K$^{-4}$. The heat transfer coefficient $q$ cannot be larger than $SL/\sigma$, as this implies a transfer of heat against the temperature gradient. \textcite{watson83} used fixed values of $q$ less than $0.2 SL/\sigma$. In this study we vary $L$ rapidly within the range of about $0.6-1.6$; for simplicity, we set $q=0.1S/\sigma$.\

As in the original model\cite{watson83} the unit of time is left unspecified. (This causes no inconsistencies with the usage of W m$^{-2}$ for fluxes due to how the equations are framed.) It is worth noting that the death rate $\gamma$ defines a natural biotic turnover timescale (i.e. timescale of biosphere renewal) $\tau_b = 1/\gamma\approx 3.3$. 

The equilibrium states of the model have been extensively studied, both numerically and analytically \cite{watson83,degregorio92,saunders94,wood08}. Nevertheless, making extensive analytical progress is challenging. One common simplification is to linearise the left-hand-side of Eq. \ref{eq:average_heat_balance}\cite{watson83,degregorio92,saunders94,wood08}. Here, we choose to focus on the original model formulation, conducting our analysis primarily using numerical continuation. To do so, we use the Continuation Core and Toolboxes (COCO) \cite{dankowicz2013recipes}. 

\begin{figure*}
	\includegraphics[width=\linewidth]{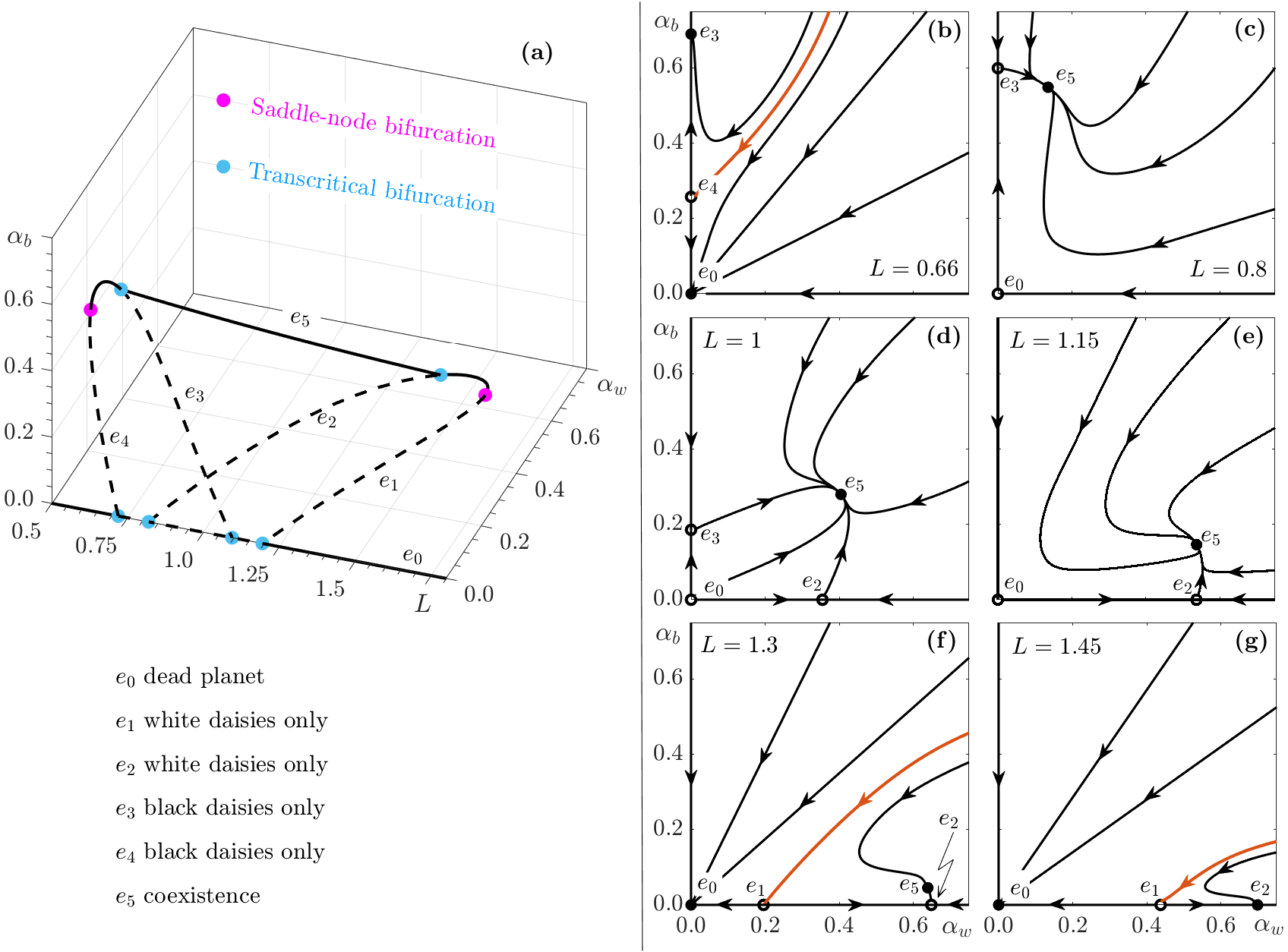}
	\caption{\label{fig:3Dbif_and_phase} \textbf{(a)} Equilibrium states in the Daisyworld model: $\alpha_b$ is the fraction of the planet covered by black daisies, $\alpha_w$ is the fraction of the planet covered by white daisies, and $L$ is the dimensionless luminosity (incoming sunlight). There are six physically relevant equilibrium states, which we label as $e_0$ through to $e_5$. $e_0$ is a dead planet, $e_1$ and $e_2$ are states in which there are only white daisies,  $e_3$ and $e_4$ are states in which there are only black daisies, and $e_5$ is a "coexistence" state, in which there are both black and white daisies. Solid lines indicate stable equilibria, and dashed lines unstable equilibria. \textbf{(b-g)} Phase portraits for different values of the luminosity $L$. Filled circles are stable equilibria, and unfilled circles are unstable equilibria. Of particular interests are the stable manifolds of certain saddle points: these are highlighted in orange 
		for easy future reference.}
\end{figure*}

Figure \ref{fig:3Dbif_and_phase} shows the equilibria of the model together with phase portraits for varying  (dimensionless) luminosity $L$. There are six equilibrium states of relevance, which we label $e_0$ through to $e_5$. $e_0$ is a dead planet (no daisies are alive at all), $e_1$ and $e_2$ are states in which there are only white daisies,  $e_3$ and $e_4$ are states in which there are only black daisies, and $e_5$ is a "coexistence" state, in which there are both black and white daisies. For a wide range of $L$ values, the coexistence state ($e_5$) is the only stable "living" state. 

The phase portraits in Figure \ref{fig:3Dbif_and_phase} further show some example phase space trajectories. Of particular interest are the stable manifolds of certain saddle points ($e_4$ for smaller values of $L$, and $e_1$ for larger values of $L$): these are highlighted in orange for easy reference later on.

\begin{figure}
\includegraphics[width=1\columnwidth]{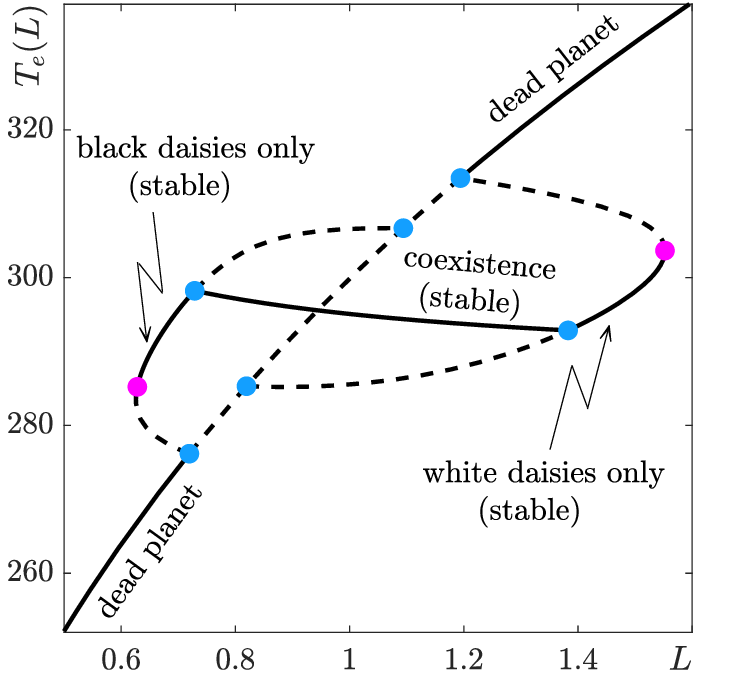}
\caption{\label{fig:bif_Te} Effective emission temperature $T_e$ as a function of varying luminosity $L$ for the different equilibrium states of the model. Solid lines are stable equilibria, dashed lines are unstable equilibria, and filled circles are bifurcations. Notably, when daisies are present, $T_e$ tends to be maintained at values that are much more conducive to the survival of daisies. For example, for larger values of $L$ (ca. 1.2-1.5) white daisies increasingly dominate, increasing the albedo of the planet and reducing overall temperature relative to a dead planet. Conversely, for smaller values of $L$ black daisies increasingly dominate, decreasing albedo and increasing overall temperature relative to a dead planet. This is the widely-discussed "self-regulation" of Daisyworld's biosphere\cite{watson83,wood08}.
}
\end{figure}

Daisyworld's biosphere is said to "self-regulate" in that, when daisies are present, the planet's effective emission temperature $T_e$ tends to be maintained at values which are much more conducive to the survival of daisies. Effectively, the daisies "make" the planet more habitable for themselves. This is demonstrated in Figure~\ref{fig:bif_Te}, which shows $T_e$ as a function of luminosity $L$ for the different equilibrium states. For larger values of $L$ the stable "coexistence" and "white daisies only" states maintain much cooler temperatures than the "dead planet" state. For smaller values of $L$ the stable "coexistence" and "black daisies only" states reliably maintain much warmer temperatures than the "dead planet" state. 

This self-regulation essentially occurs due to the interplay of three factors: the daisies' temperature-dependent growth rates; the daisies' differing albedos; and the fact that, due to these albedo differences, black daisies will be warmer than white daisies. If $L$ is relatively small, the black daisies will outcompete the white daisies, tending to reduce the albedo relative to uncovered ground and warming the planet. Conversely, if $L$ is relatively large, the white daisies will outcompete the black daisies, tending to increase the albedo relative to uncovered ground and cooling the planet.

\section{Bifurcation-induced collapse}
We now consider tipping points in the Daisyworld model. In particular, we focus on changes in the dimensionless luminosity $L$ which may cause a "living" equilibrium to abruptly go extinct. The first kind of tipping we consider is bifurcation-induced tipping\cite{ashwin12}. 

Conceptually, bifurcation-induced tipping can be understood by considering the general nonautonomous system 
\begin{equation}
	\di{x}{t} = f(x,\Lambda(rt)),
	\label{eq:open_system}
\end{equation}

in which $\Lambda(r t)$ is an input parameter that varies on a much slower time scale $r\,t$. Here, $r$ 
is a small parameter whose unit is 1/(unit time).
In this setting, the parameter $r$ quantifies the rate of change of the input parameter $\Lambda$, and is often referred to as
the {\em rate parameter} \cite{ashwin2017parameter, ritchie23}.

If the autonomous system 
\begin{equation}
	\di{x}{t} = f(x,\lambda)
\end{equation}
exhibits a "dangerous" bifurcation\cite{thompson1994safe} at some $\lambda=\lambda^*$, then the nonautonomous system in Eq. \ref{eq:open_system} may tip when $\Lambda(r t)$ crosses $\lambda^*$. Specifically, this will occur if the system state is initially in a stable equilibrium which loses stability as a result of the bifurcation, and $\Lambda(r t)$ remains beyond $\lambda^*$ for a sufficiently long time interval. As discussed in the Introduction, bifurcation-induced tipping is responsible for fixed (i.e. time-independent) thresholds ($\lambda^*$in this case). In this case, as long as $\Lambda(rt)$ crosses $\lambda^*$, the precise value of $r$ does not affect whether or not tipping occurs.

\begin{figure*}
	\includegraphics[width=0.9\linewidth]{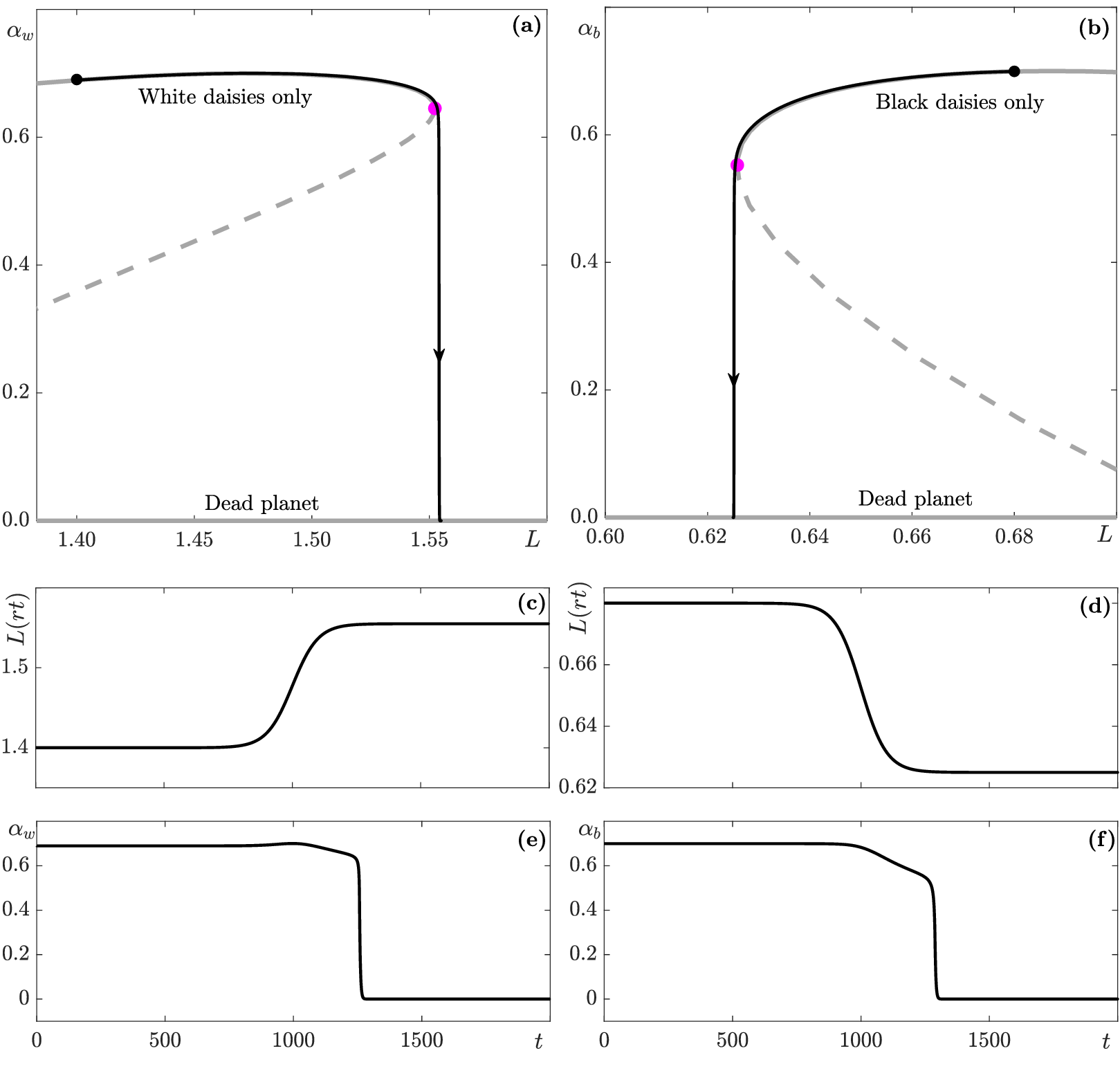}
	\caption{\label{fig:2Dbif_btipping} Bifurcation-induced tipping towards extinction in the Daisyworld model. In both (a) and (b) we initialise the system at a stable "living" equilibrium (black dot), and then gradually change the luminosity $L(rt)$. In panel (a), there is initially a large population of white daisies; however, eventually the equilibrium undergoes a saddle-node bifurcation (magenta dot), the only remaining stable equilibrium is that of a dead planet, and the biosphere collapses. In panel (b) we demonstrate the same phenomenon starting from the stable state in which there are only black daisies, and decreasing $L$. In each case, dashed lines denote unstable equilibrium states.
    Panels (c), (d) show time series of the luminosity $L(rt)$ with $r = 0.01$, corresponding to (a) and (b), respectively. Panels (e) and (f) show time series of $\alpha_w$ and $\alpha_b$ corresponding to (a) and (b), respectively.}
\end{figure*}

In the Daisyworld model, increasing $L$ to a large enough value, or decreasing it to a small enough value, will lead to bifurcation-induced tipping towards extinction. This is already apparent from the left-hand panel in Figure \ref{fig:3Dbif_and_phase}, and we demonstrate this in more detail in Figure \ref{fig:2Dbif_btipping}. Here and throughout the rest of the paper, we implement time-varying forcings in luminosity $L$ according to
\begin{equation}
	L(rt) = L_{\rm start} + \frac{\Delta L}{2} (\tanh(r t) + 1),
	\label{eq:Lt_function}
\end{equation}
where $\Delta L=L_{\rm end}-L_{\rm start}$, and the rate parameter $r$ controls how quickly the transition between $L_{\rm end}$ and $L_{\rm start}$ occurs. Unless specified otherwise, we will use values of $L_{\rm end}>L_{\rm start}$.

In Figure \ref{fig:2Dbif_btipping} we consider two scenarios. In the first scenario, shown in Figure~\ref{fig:2Dbif_btipping} (a), (c), and (e), we initialise the system at fairly high luminosity $(L=1.4)$, from the stable equilibrium in which only white daisies exist. Then, as luminosity is increased, eventually the stable equilibrium undergoes a saddle-node bifurcation, and only the "dead planet" state is now stable. As a result, the system inevitably evolves towards the dead planet state, and the daisies go extinct. Physically, the planet's surface temperature was initially maintained by the white daisies far below where it would otherwise have been. However, the effect of increased temperature on daisy growth rates (Eq. \ref{eq:growth_rates}) eventually overwhelms the stabilising effect of the daisy population on temperature, and the daisy population abruptly collapses.

In the second scenario, shown in Figure~\ref{fig:2Dbif_btipping} (b), (d), and (f), we initialise the system at fairly low luminosity $(L=0.68)$, from a state in which only black daisies exist, We now decrease the luminosity (in this particular case $L_{\rm end}<L_{\rm start}$) until the equilibrium disappears through a saddle node bifurcation, with similar consequences as in the first case. Initially the planet's temperature was maintained by the black daisies far above where it would otherwise have been. But again, the effect of temperature on daisy growth rates eventually overwhelms the stabilising effect of the daisy population on temperature, and the daisy population collapses.

There is no bifurcation-induced tipping starting directly from the coexistence equilibrium state, because the only bifurcations this can undergo are transcritical bifurcations to the single-species equilibria (see Figure \ref{fig:3Dbif_and_phase}). However, we note that if the system is initialised at the coexistence equilibrium state for any given value of $L$, and $L$ is slowly increased or decreased, for $|\Delta L|$ large enough the system will eventually undergo bifurcation-induced tipping: it is just that the system will pass through a single-species equilibrium first. 

\section{Rate-induced collapse and basin instability}
\begin{figure*}
	\includegraphics[width=\linewidth]{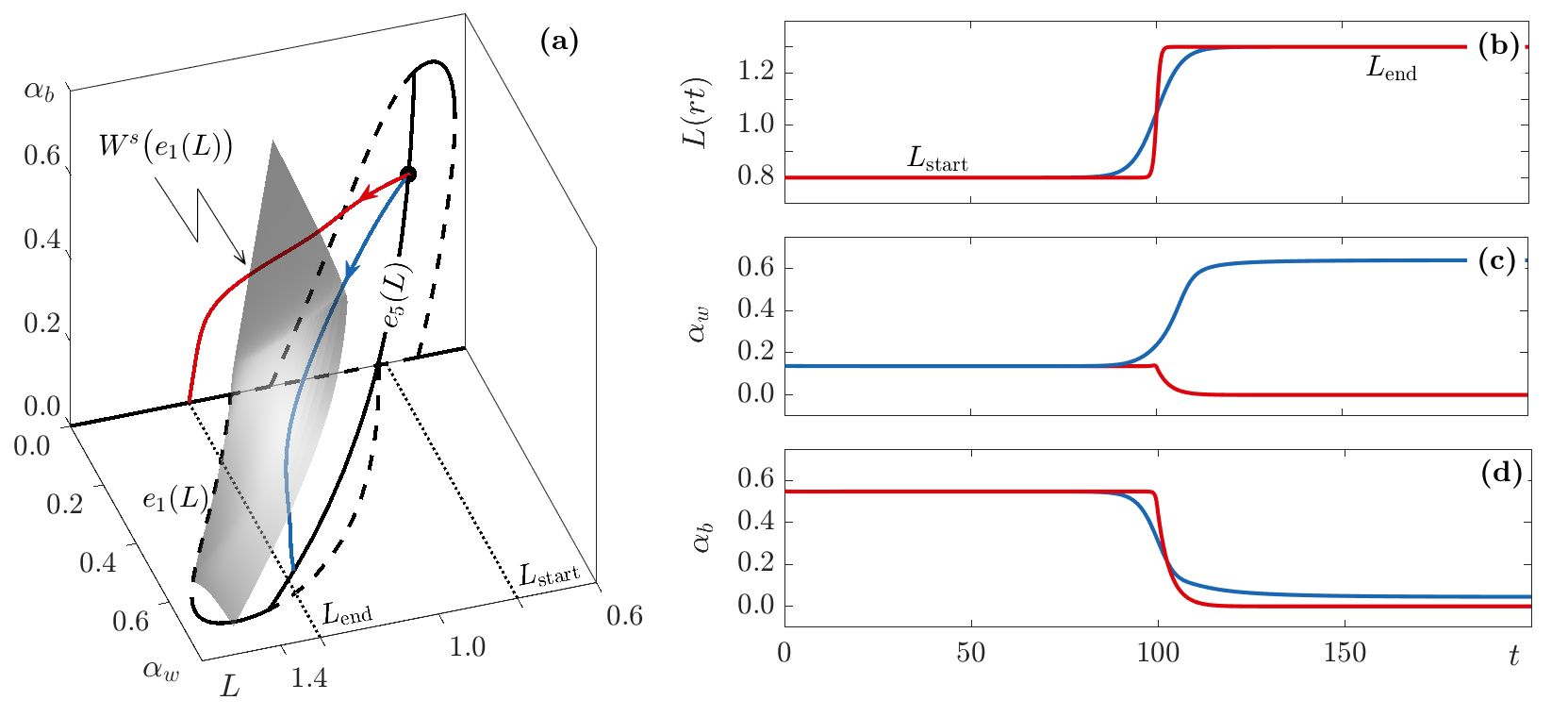}
	\caption{\label{fig:3Dbif_rtipping} Rate-induced tipping towards extinction in the Daisyworld model. We initialise the system at the coexistence equilibrium $(e_5)$ with a value of $L=L_{\rm start}$. Then, we ramp $L$ towards some new value $L_{\rm end}$, at a speed determined by a dimensionless rate parameter $r$. Panel (b) shows time series of $L(rt)$, while panels (c) and (d) show time series of $\alpha_w$ and $\alpha_b$. When the change is slow enough (blue, $r = 0.5$), the biosphere survives; when it is too fast (red, $r = 1$), the biosphere collapses. This occurs because the system can enter the basin of attraction of the "dead planet" state, $e_0$, if the change in $L$ is fast enough. The threshold for entering this basin of attraction is the stable manifold of the saddle equilibrium $e_1$ (see also Figure~\ref{fig:BI_phase} and panel f of Figure \ref{fig:3Dbif_and_phase}). In the space of $\alpha_w$,$\alpha_b$, and $L$, the threshold becomes a two dimensional surface $W^{s}(e_1(L))$: trajectories which cross it are precisely those which tip. This is shown in panel (a), in which the red curve crosses the surface and the blue curve does not.}
\end{figure*}

Next, and more notably, we show that Daisyworld's biosphere can undergo rate-induced tipping towards extinction.  This is demonstrated in Figure \ref{fig:3Dbif_rtipping}.  We initialise the system at the coexistence equilibrium with some smaller $L$ value $L_{\rm start}$, and consider 
the time-dependent forcings $L(rt)$ (Eq. \ref{eq:Lt_function}) which ultimately increase $L$ to some larger value $L_{\rm end}$. The coexistence equilibrium exists and is stable in the entire range $L_{\rm start}\leq L\leq L_{\rm end}$. 

In Figure~\ref{fig:3Dbif_rtipping} we show the result of such a perturbation for two different values of $r$: one fast (red), and one slow (blue). In both cases, the change is fast enough that the system state fails to track the moving equilibrium ($e_5(L)$). In the slow case, the system is able to eventually reach this equilibrium ($e_5(L_{\rm end})$): the biosphere has survived the change in $L$. However, when $L$ changes more quickly (red), the biosphere collapses, leaving a dead planet. The collapse occurs even though a stable biosphere could in principle exist across the entire range $L_{\rm start}\leq L\leq L_{\rm end}$. (i.e. the coexistence equilibrium is stable for all $L_{\rm start}\leq L\leq L_{\rm end}$). This is rate-induced tipping.

\begin{figure}
\includegraphics[width=0.95\columnwidth]{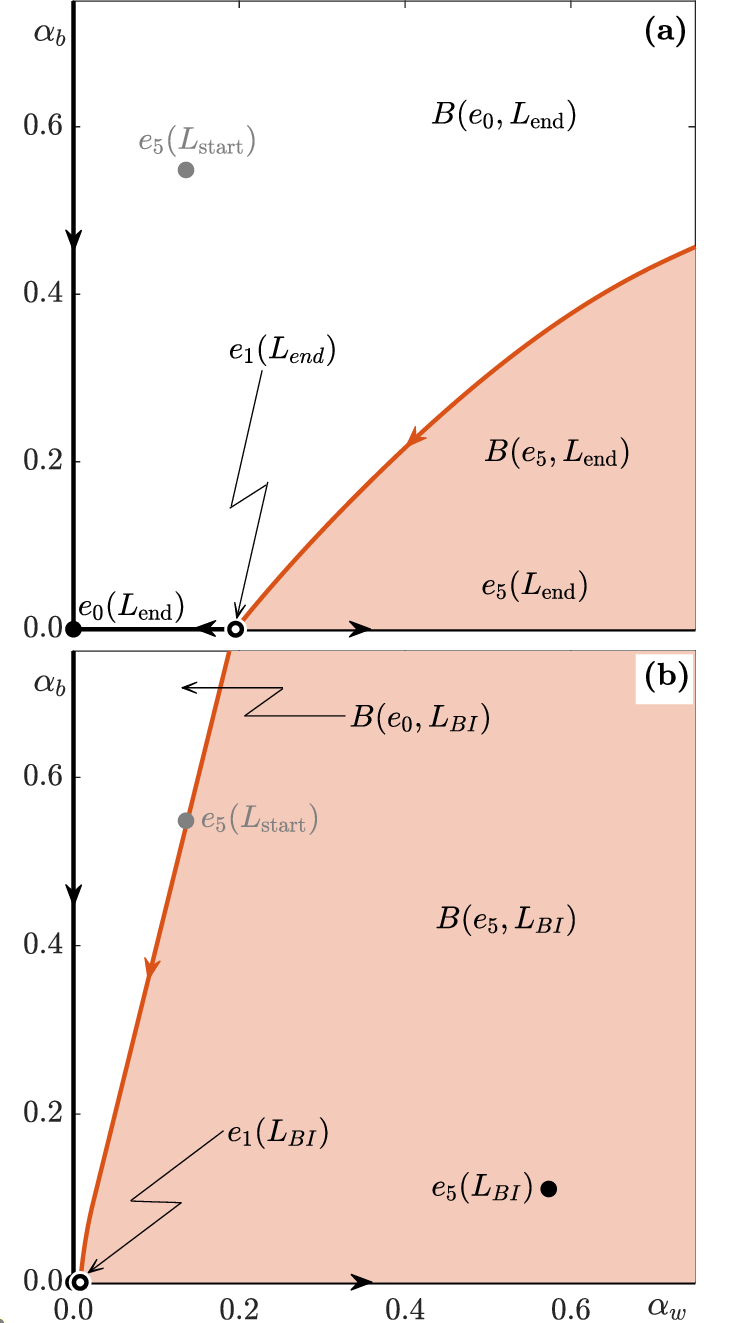}
\caption{\label{fig:BI_phase} Rate-induced tipping in the Daisyworld model (e.g. as seen in Figure \ref{fig:3Dbif_rtipping}) is a consequence of basin instability. Panel (a) shows that, at $L=L_{\rm end}$ (the upper end of the perturbation in Figure \ref{fig:3Dbif_rtipping}) the phase space is partitioned into two basins of attraction, separated by the stable manifold of the saddle point $e_1(L_{\rm end})$. The stable coexistence equilibrium $e_5(L_{\rm start})$ lies within the basin of attraction of $e_0(L_{\rm end})$. Thus, if the system is initialised at $e_5(L_{\rm start})$ with $L=L_{\rm start}$, and then $L$ is instantaneously changed to $L_{\rm end}$, the system tips towards $e_0$ (the "dead planet" state). The equilibrium $e_5$ is thus basin unstable\cite{okeeffe20}. Panel (b) shows the two basins at $L=L_{BI}$, which is where $e_5(L_{\rm start})$ intersects the boundary of the basins. We refer to this as the boundary of the basin instability region of equilibrium $e_5(L_{\rm start})$. 
}

\end{figure}

To understand rate-induced tipping, we consider Figure~\ref{fig:BI_phase}, which shows the coexistence stable equilibrium $e_5(L)$, the dead planet equilibrium $e_0(L)$, and their basins of attraction at different values of the parameter $L$, for the autonomous Daisyworld model. 

As mentioned before, the coexistence equilibrium $e_5(L)$ is stable for luminosity $L$ in $[L_{\rm start}, L_{\rm end}]$, and so is the dead planet equilibrium $e_0(L)$.
This means that for all $L$ in $[L_{\rm start}, L_{\rm end}]$ the phase space is partitioned into two basins of attraction. These are the basins of attraction of $e_5(L)$ and of $e_0(L)$.
To simplify the discussion, we write
$$
B(e_5,L),
$$
to denote the basin of attraction of $e_5$ at the parameter value $L$, and 
$$
B(e_0,L),
$$
to denote the basin of attraction of $e_0$ at the parameter value $L$. 
The boundary between these two basins is the stable manifold of the saddle equilibrium $e_1(L)$ (orange curve). 

In Figure~\ref{fig:BI_phase} we plot $B(e_5,L)$ (light orange) and $B(e_0,L)$ (white) for two different luminosity values $L=L_{\rm end}$ and $L=L_{BI}$, Additionally, we show the stable coexistence equilibrium $e_5(L_{\rm start})$ as a grey point, which represents the starting point for the R-tipping trajectories in Figure~\ref{fig:3Dbif_rtipping}. 

In panel (a) for $L=L_{\rm end}$, the equilibrium $e_5(L_{\rm start})$ lies outside the basin of attraction of $e_5(L_{\rm end})$, i.e. 
$$
e_5(L_{\rm start}) \not\in B(e_5,L_{\rm end}). 
$$
In such a case, the equilibrium $e_5(L_{\rm start})$ is said to be {\em basin unstable} on the parameter path $[L_{\rm start}, L_{\rm end}]$ \cite{okeeffe20}. 

Figure~\ref{fig:BI_phase}~(b) shows the two basins at $L = L_{BI}$, which is the value of the parameter where the equilibrium $e_5(L_{\rm start})$ intersects the boundary of the basins. We refer to $L_{BI}$ as the boundary of the basin instability region of equilibrium $e_5(L_{\rm start})$.  

To explain rate-induced tipping, when $L(rt)$ changes from $L_{\rm start}$ to $L_{\rm end}$, we consider two extreme cases. First, when $r \to 0$, $L( r t)$ changes extremely slowly. Thus, the system can keep track of the moving equilibrium $e_5(L(r t))$. 
Second, when $r \to \infty$, $L( r t)$ changes instantaneously from $L_{\rm start}$ to $L_{\rm end}$. Hence, the state of the system, although it starts at $e_5(L_{\rm start})$, ends up in the basin of $e_0(L_{\rm end})$. Therefore, the system tips.
By continuity, there is a critical value of the rate parameter $r_c$ such that the system exhibits tipping for any $r > r_c$ and tracks for $r < r_c$. 
If $r = r_c$ the system converges to the saddle point $e_1(L_{\rm end})$\cite{wieczorek2023rate}. (We note that transition behaviour in dynamical systems due to crossing stable manifolds of saddle points has long been considered, see e.g. \textcite{fitzhugh55}). 

One can consider an augmented (compactified) system where $L(rt)$ can be incorporated as an additional differential equation. 
Such an augmented system will exhibit a codimension-one heteroclinic bifurcation for $r = r_c$ \cite{wieczorek2023rate}. 
For sufficiently simple models the heteroclinic orbit, and hence $r_c$, can be computed in closed form \cite{ritchie2016early} or approximated using perturbation methods. 
To highlight rate-induced biosphere collapse in the Daisyworld model, we find it sufficient to provide a numerical approximation of the critical rate parameter $r_c$.

When $L$ changes rapidly ($r > r_c$), the system state moves from the basin of attraction of the coexistence state ($e_5$) to the basin of attraction of the dead planet state ($e_0$). The threshold boundary --- the stable manifold of $e_1$ --- is shown in Figure \ref{fig:3Dbif_rtipping}: in the space of $\alpha_w,\alpha_b$ and $L$, it becomes a two-dimensional surface. Those trajectories which undergo rate-induced tipping are precisely those which pass through this surface.

\begin{figure*}
\includegraphics[width=0.95\linewidth]{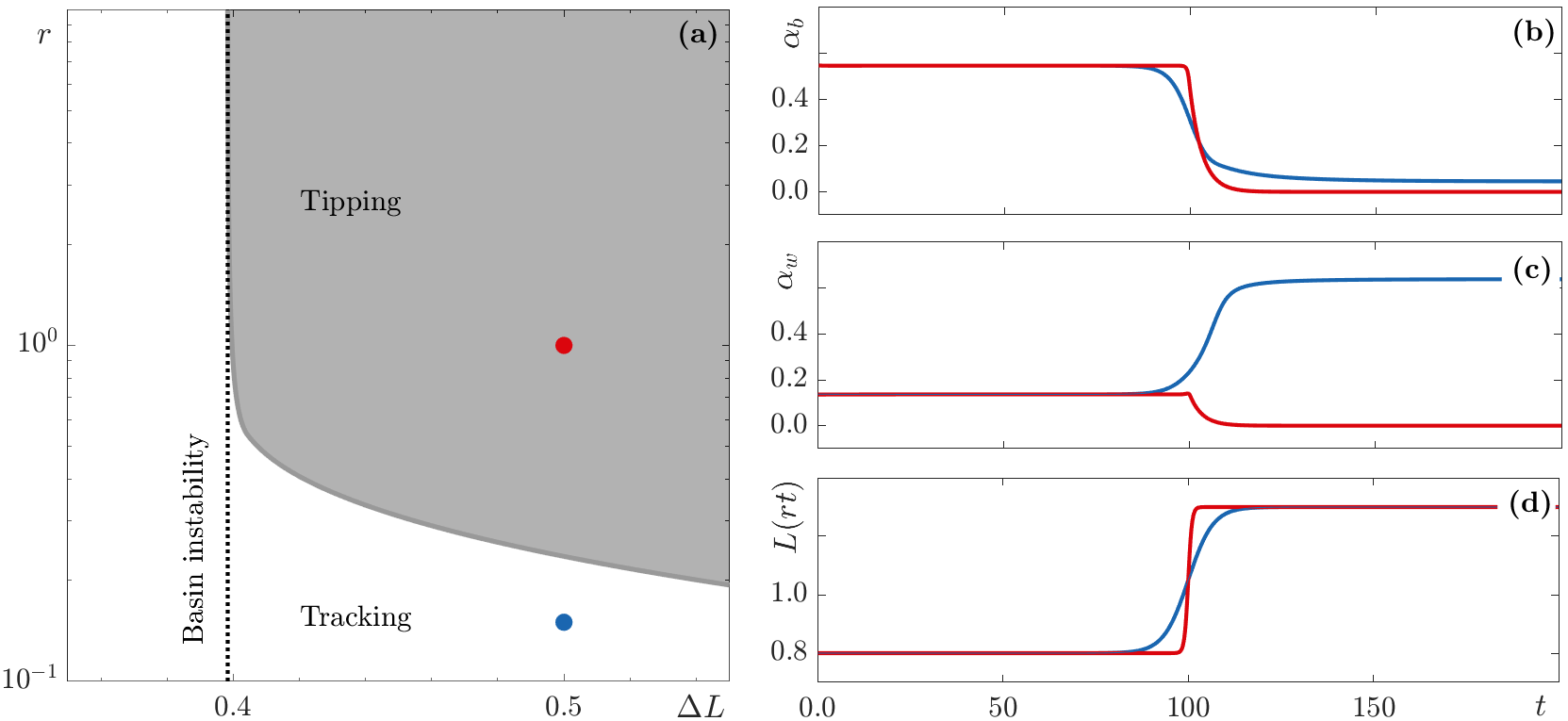}
\caption{\label{fig:Rtipping_diagrams}Tipping diagram for the Daisyworld model, starting from the same coexistence equilibrium $e_5$ with luminosity $L=L_{\rm start}$. Panels (b), (c), and (d) again show time series of $\alpha_b$,$\alpha_w$, and $L(rt)$, for the same two examples as in Figure \ref{fig:3Dbif_rtipping}. Panel (a) shows, for different combinations of rate $r$ and total luminosity change $\Delta L = L_{\rm end} - L_{\rm start}$, which perturbations result in tipping, and which ones do not. For large rates $r$, the critical $\Delta L$ converges to a fixed value, which represents the instantaneous $\Delta L$ required to move the system across the basin threshold.
}
\end{figure*}

We can further characterise the kinds of perturbations that trigger rate-induced tipping. We continue to initialise the system at the stable coexistence equilibrium $e_5(L_{\rm start})$, and perturb it using the forcing function in Eq. \ref{eq:Lt_function}, varying both the rate parameter $r$ and the total change in $L$, $\Delta L = L_{\rm end} - L_{\rm start}$.
Our results are shown in Figure \ref{fig:Rtipping_diagrams}, with the two examples from Figure \ref{fig:3Dbif_rtipping} again highlighted. For fast enough rates $r$, the  $\Delta L$ required to tip the system converges to a constant value, which represents the $\Delta L$ value required to instantaneously move the system state across the basin instability boundary. For slower perturbations, the system is able to remain in the original basin of attraction, and is thus always able to recover to the original (but now shifted) coexistence state. 

One further conclusion from this diagram is that, on Daisyworld, for a general perturbation in $L$, basin instability (i.e. rate-induced tipping) is often "easier" to initiate than bifurcation-induced tipping. To see this, we begin by looking at Figure \ref{fig:3Dbif_rtipping}, and noting that we start from a value of $L_{\rm start} = 0.8$. For an increase in $L$ which is monotonic but whose time dependence is otherwise arbitrary, the furthest that $L$ could possibly increase while Daisyworld's biosphere survives is until the saddle-node bifurcation at $L\simeq1.55$: see Fig. \ref{fig:3Dbif_and_phase} panel a, and Fig \ref{fig:2Dbif_btipping}, panel a. (We note that, by this point, the stable living equilibrium is no longer the coexistence state $e_5$, but instead the "white daisies only" state $e_2$.) Thus, the point at which $\Delta L$ must necessarily trigger extinction is at $\Delta L \simeq 1.55-0.8=0.75$. But, as Figure \ref{fig:3Dbif_rtipping} shows, perturbations of sufficient speed can trigger extinction for much smaller values of $\Delta L$. On Daisyworld, rates of change can play as important a role in determining catastrophic outcomes as the magnitude of the change. 

\section{Discussion}
We have shown that the Daisyworld model exhibits both bifurcation-induced tipping to extinction, and rate-induced tipping to extinction via basin instability. We emphasise that the precise relationship of the Daisyworld model to Earth's real-world biosphere is a contentious subject \cite{kirchner89,kirchner02,wood08}. Regardless, our result is of interest, for a number of reasons.

First, it improves the paedagogical value of the Daisyworld model. Daisyworld has already been widely used in the teaching of Earth system science and co-evolution\cite{mcguffie05_book,kump10_book,ford10_book}, and could now be used to introduce students to an even richer range of dynamics, like rate-induced tipping, which are also relevant for understanding various real-world systems.

Second, this result serves as another useful data point in the ongoing discovery of rate-induced tipping in more and more systems. Recent years have seen an explosion of evidence regarding the ubiquity of rate-induced tipping, both in models and in the real world \cite{stocker1997influence, scheffer08,wieczorek10, luke2011soil, siteur2016ecosystems, alkhayuon2019basin, okeeffe20, vanselow19,lohmann21,arnscheidt20,arnscheidt22_review,arnscheidt22_rate,ritchie23, feudel2023rate, panahi2023rate}. If a classic four-decade old model like the Daisyworld model can, upon re-assessment, be shown to exhibit rate-induced tipping, which other classic models might this hold for?

Third, this helps support the argument that rate-induced collapse is common in systems with some capacity to evolve and/or adapt\cite{arnscheidt22_rate,arnscheidt22_review}. It is well understood that such adaptive capacity may allow a system to survive under slowly changing external conditions. But if this is the case, the same system can fail to survive when change is too fast: in other words, it undergoes rate-induced tipping towards collapse. This is exactly the behaviour we see on Daisyworld, and which we have described in this article. We thus hope that, beyond the specifics,  this work can serve as yet another useful example of a phenomenon which likely has much broader relevance, especially as we continue to navigate a world dominated by rapid human-driven rates of change.  

\section*{Code Availability}
The codes used to conduct simulations and generate ﬁgures are available via the GitHub repository: \\
\url{https://github.com/hassanalkhayuon/Daisyworld_}

\begin{acknowledgments}
Both authors thank the organisers of the October 2023 workshop "Non-autonomous Dynamics in Complex
Systems: Theory and Applications to
Critical Transitions" for helping catalyse this collaboration. C. W. A. thanks D. Rothman, G. Gibbins, and E. Stansifer for helpful discussions on an earlier version of this work. 
\end{acknowledgments}

\bibliography{daisyworld_2}

\end{document}